\documentclass[12pt]{article}
\usepackage{amsmath}
\usepackage{amssymb, amsbsy, amstext}
\usepackage{srcltx}
\usepackage{amscd,amsfonts,color}
\usepackage{float}

\usepackage{graphicx}
\usepackage{stmaryrd}
\usepackage{mathabx}
\usepackage{enumerate}
\usepackage{url}
\usepackage{tikz}
\usepackage{pgflibraryarrows}
\usepackage{pgflibrarysnakes}

\input amssym.def
\input amssym.tex

\newtheorem{lem}{Lemma}[section]%
\newtheorem{theo}[lem]{Theorem}%
\newtheorem{defi}[lem]{Definition}%
\newtheorem{cor}[lem]{Corollary}%
\newtheorem{prob}[lem]{Problem}%
\newtheorem{prop}[lem]{Proposition}%

\newcommand{\ZZ}{\mathbb{Z}}

\newcommand{\B}{{\rm B}}

 \newcommand{\ox}{\overline X}

\newcommand{\Aut}{\hbox{\rm Aut\,}}

\newcommand{\Core}{\hbox{\rm Core}}

\newcommand{\demo}{\f {\bf Proof}\hskip10pt}

\newcommand{\BB}{\mathbf{B}}
\newcommand{\Bc}{\mathcal{B}}

\newcommand{\ux}{\underline{X}}

\newcommand{\A}{{\mathcal{A}}}
\newcommand{\ba}{{\bf \alpha}}

\def\a{\alpha}   \def\d{\delta} 
 \def\s{\sigma}

\def\si{\Sigma}    
\def\ol1{\overline 1}
\def\di{\bigm|}
\def\lg{\langle}
\def\rg{\rangle}

\def\f{\noindent}
\def\qed{\hfill $\Box$}

\begin{document}
\begin{center} {\bf\large Hamilton Cycles in Semisymmetric Graphs}
\bigskip

{\sc Shaofei  Du$^\dag$ and Kai Yuan$^*$} \\

{\footnotesize  $^\dag$School of Mathematical Sciences,  Capital Normal University, Beijing 100048, China}\\
{\footnotesize   $^*$School of Mathematics and Information Sciences, Yantai University, Yantai 264005, China }\\

\end{center}

\begin{abstract} In light of Lov\'{a}sz's longstanding question on the existence of Hamilton paths in vertex-transitive graphs, this paper considers a natural variant: what if vertex-transitivity is relaxed, yet a high degree of symmetry--specifically edge-transitivity--is retained? To investigate this, we focus on the class of semisymmetric graphs, which are regular, edge-transitive, but not vertex-transitive. 
In this paper, it will be  shown   that  every  connected semisymmetric graph of order $2pq$, where $p$ and $q$ are two distinct  primes contains a Hamilton cycle  and that  every connected cubic semisymmetric  graph of order less than 3000  
contains a Hamilton cycle too. Based on these observations,  the following question is posed:
construct a  connected semisymmetric graph which has no Hamilton cycle.
   \end{abstract}

{\small  {\em 2000 Mathematics Subject Classification}:\,
             05C25, 05C45.

             }

\renewcommand{\thefootnote}{\empty}
\footnotetext{Corresponding author: pktide@163.com (Kai Yuan).}
\footnotetext{E-mail addresses: dushf@mail.cnu.edu.cn (Shaofei Du), pktide@163.com (Kai Yuan).}
\footnote{This work was supported by the National Natural Science Foundation of China
 (12471332 and 12101535).}
 \section{Introduction}
The graphs considered in this paper are finite, simple, undirected. A {\it Hamilton	path $($resp. cycle$)$}	of a graph  is  a simple path (resp. cycle) containing	 all 	vertices of the graph.
The following   question asked by Lov\'{a}sz \cite{Lo} in 1969 tying together traversability and symmetry
remains unresolved after all these years:

\vskip 3mm  {\it Construct a finite connected vertex-transitive graph which has no Hamilton path.}
\vskip 3mm

No	connected vertex-transitive graph without a Hamilton	path
 is known to exist and
so far  only four
connected vertex-transitive  graphs	on at least three	vertices not having  a Hamilton	cycle are	
known: Petersen graph, the Coxeter graph, and  another two respective   truncations from them.  Since  none of these four graphs is Cayley graph. 	
This has led to another conjecture that every connected Cayley graph
on at least three vertices possesses a Hamilton cycle.

 There have been
 a	number of papers affirming the existence of Hamilton 	paths and in some cases even Hamilton cycles for some  families of vertex transitive graphs.
These graphs contain:   connected vertex-transitive graphs of orders	
$kp, k \leq 6$, $p^j, j \leq 5$ and $2p^2$,
$10p, p\ge 11$, 
where	$p$ is	a prime	 (see
\cite{A79,C98, DZ, DLY, KM08, KMZ12,KS09,DM85,DM87,MP82,MP83,Z15}
and the survey paper \cite{KM09}),  except for the graphs of orders  $10p$
(and of course for  the above four graphs),
it is also known that they contain a Hamilton cycle;  and connected
vertex-transitive graphs whose automorphism groups contain a transitive subgroup with a cyclic commutator subgroup of prime-power order (except for  the Petersen graph) \cite{DGMW98};
 and some Cayley graphs with given  groups (see
\cite{A89,GKMM12,GM07,DM83,WM18} and the survey papers \cite{CG96,KM09}).
Very recently, with  long term efforts, Du,  Kutnar and Maru\v si\v c  proved  the existence of Hamilton 	cycles  of connected vertex-transitive graph of
order $pq$, where $p$ and $q$ are primes, with the exception of the Petersen graph \cite{DKM}.

 \vskip 3mm
 Motivated by  Lov\'{a}sz  Question, we may ask  the following question ?
 \vskip 3mm

 {\it How about  the Hamiltonicity for graphs,  if we remove its vertex-transitivity  but still keep a high symmetric property for it ? }
 \vskip 3mm

 To investigate such problem, we shall consider a class of graphs: {\it semisymmetric graphs}, which are regular,  edge-transitive but not vertex-transitive.

Let's recall some backgrounds of semisymmetric graphs.  It is easy to see from its definition  that every
 semisymmetric graph is a bipartite graph with two biparts of equal size  and the automorphism group acts transitively of every  bipart.
 In 1967, Folkman firstly constructed several infinite families of semisymmetric graphs and proposed eight open problems (see \cite{Fo}).
   In paticular, he proved that there are no semisymmetric graphs of orders $2p$ and $2p^2$.
    Afterwards, many authors have done quite much work on this topic.  Among of them,  depending on classification of finite simple groups, Iofinova and Ivanov \cite{II} in 1985 classified all  cubic semisymmetric graphs whose automorphism group acts primitively on  every bipart (so  called {\it biprimitive}) and they proved that there are only five such graphs; and   Du and Xu  \cite{DX}  classified semisymmetric graphs of order $2pq$, where $p$ and $q$ are distinct primes. Moreover,
     in \cite{CP}, Conder and Poto\v{c}nik generate all the connected semi-symmetric 3-valent graphs of order up to 10000.

  Now we are ready to state the main results of this paper.
\begin{theo}\label{a}
Every  connected semi-symmetric cubic graph of order less than 3000 has a   Hamilton cycle.
\end{theo}
\begin{cor} \label{b}
Every  cubic biprimitive semisymmetric graph has a   Hamilton cycle.
\end{cor}
\begin{theo}\label{c}
Every connected semisymmetric graph of order $2pq$ contains a Hamilton cycle, where $p$ and $q$ are distinct  primes.
\end{theo}
Based on these observations, we would like to pose a problem:
\begin{prob}\label{d}
Construct a finite connected semisymmetric graph which has no Hamilton cycle.
\end{prob}

\section{Preliminaries}
In this section, we recall some definitions and known results, while all of them are quoted from \cite{DX}.

 If $X$
is bipartite with bipartition $V(X)=U(X)\cup W(X)$, we let
$A^+=\lg g\in A\di U(X)^g=U(X), W(X)^g=W(X)\rg$.   Clearly if $X$ is
connected then either $|A:A^+|=2$ or $A=A^+$, depending on whether
or not there exists an automorphism which interchanges the two parts
$U(X)$ and $W(X)$.
Suppose $G$ is a subgroup of $A^+$. Then $X$ is said to be
$G$-{\em semitransitive} if $G$ acts transitively on both $U(X)$
and $W(X)$, and {\em semitransitive} if $X$ is $A^+$-semitransitive.
Also $X$ is said to be {\it biprimitive} if $A^+$ acts primitively
on both $U(X)$ and $W(X)$.

\begin{defi}
Let $G$ be a group,   let $L$ and $R$ be   subgroups of $G$
and let  $D$ be a union of double cosets of $R$ and $L$ in $G$, namely,
$D=\bigcup_i Rd_iL$.
Define a bipartite graph $X=\B(G,L,R;D)$ with bipartition $V(X)=[G:L]\cup
[G:R]$ and  edge set $E(X)=\{\{Lg,Rdg\}\di g\in G, d\in D\}$.
This graph is called the bi-coset graph of $G$ with respect to $L$, $R$ and
$D$.
\end{defi}

\begin{prop}\cite{DX} \label{property}
The graph $X=\B(G,L,R;D)$ is a well-defined bipartite graph.
Under the right multiplication action on $V(X)$ of $G$, the graph $X$ is
$G$-semitransitive.  The kernel of the action of $G$ on $V(X)$ is
$\Core_G(L)\cap\Core_G(R)$, the intersection of the cores of the
subgroups $L$ and $R$ in $G$. Furthermore, we have

\item{\rm (i)} $X$ is $G$-edge-transitive if and only if $D=RdL$
for some $d\in G$;

\item{\rm (ii)} the degree of any vertex in $[G:L]$ (resp. $[G:R])$
is equal to the number of right cosets of $R$ (resp. $L$) in $D$
(resp. $D^{-1})$, so $X$ is regular if and only if $|L|=|R|$;

\item{\rm (iii)} $X$ is connected if and only if $G$ is generated
by elements of $D^{-1}D$;

\item{\rm (iv)} $X\cong \B(G,L^a,R^b;D')$
where $D'=\bigcup_iR^b$ $(b^{-1}d_ia)L^a$, for any $a, b\in G$.
\end{prop}

\begin{prop}\cite{DX} \label{property2}
Suppose $Y$ is a $G$-semitransitive graph with
bipartition $V(Y)=U(Y)$ $\cup\, W(Y)$.
Take $u\in U(Y)$ and $w\in W(Y)$. Set \
$D=\{ g\in G\di w^g \in Y_1(u)\}.$ \
Then $D$ is a union of double cosets of $G_w$ and
$G_u$ in $G$, and $Y\cong \B(G,G_u,$ $G_w;D).$
\end{prop}

By Propositions \ref{property} and \ref{property2}, every $G$-semitransitive and  $G$-edge-transitive graph  can be represented by a bi-coset graph  $\B(G,G_u,$ $G_w;D),$ where 
  $u\in U(Y)$,  $w\in W(Y)$ and $D=G_wdG_u$ for an edge $\{u, w^d\}$.

\begin{defi}

Given a $G$-semitransitive and $G$-edge-transitive graph $X$ with bipartition $V(X)=U(X)\cup W(X)$, where $|U(X)|=p$
and $|W(X)|=pq$ for two distinct prims $p$ and $q$, we define a bipartite graph $\ux$ with bipartition $U(\ux)\cup W(\ux)$
as follows:
$$\begin{array}{lll} U(\ux)&=&\{u_i \di u \in U(X), 1 \leq i \leq q \},~~~~ W(\ux)=W(X),\\
E(\ux)&=&\{(u_i,w) \di (u,w)\in E(X), 1 \leq i \leq q\}.\end{array} $$
We call $\ux$ the derived graph of $X$.
\end{defi}

The following proposition gives a classification of semisymmetric graphs of order $2pq$.
\begin{prop} \cite[Theorem 2.8]{DX}
Suppose $Y$ is a semisymmetric graph of order $2pq$, where $p$ and $q$ are primes such that $q<p$. Then the following hold:
\begin{enumerate}
\item [$(i)$] If $Y$ is biprimitive, then $Y$ is isomorphic to one of the bi-coset graphs $\BB(G,L,R;$ $D)$ given in TABLE 1 in terms of
the triple $(G,L,R)$, the primes $p$ and $q$, the valency of $Y$, and the number of nonisomorphic possibilities with the same parameters.
\item [$(ii)$] Otherwise, if $Y$ is not biprimitive, then $Y$ is isomorphic to the graph $\ux$ derived from one of the bi-coset graphs
$X=\BB(G,L,R;D)$ given in TABLE 2 in terms of $(G,L,R)$, $p$ and $q$, the degree of a vertex $u \in [G:L]$, and the number of nonisomorphic possibilities with the same parameters.
\end{enumerate}
\end{prop}

\begin{center}
TABLE 1: Biprimitive bi-coset graphs $\BB(G,L,R;$ $D)$
\vskip 5mm
{\small
\begin{tabular}{|c|c|c|c|c|} \hline
Row &$(G, L, R)$ &$pq$ & Valency & Number \\ \hline \hline
1 &$(PSL(2, 59), A_5, A_5)$ &$59\cdot 29$ &60 &4 \\ \hline
2 &$(PSL(2, 61), A_5, A_5)$ &$61\cdot 31$ &60 &4 \\ \hline
3 &$(PSL(2, 23), S_4, S_4)$ &$23\cdot 11$ &12, 24 &1, 1 \\ \hline
4 &$(PSL(2, 59), A_5, A_5^{\s })$ &$59\cdot 29$ &60 &4 \\ \hline
5 &$(PSL(2, 61), A_5, A_5^{\s })$ &$61\cdot 31$ &60 &5 \\ \hline
6 &$(PSL(2, 23), S_4, S_4^{\s })$ &$23\cdot 11$ &12, 24 &1, 1 \\ \hline
7 &$(PGL(2, 11), S_4, D_{24})$ &$11\cdot 5$ &3, 4, 12, 24 &1, 1, 2, 1\\ \hline
8 &$(PSL(2, 13), A_4, D_{12})$ &$13\cdot 7$ &6, 12 &1, 2 \\ \hline
9 &$(PGL(2, 13), S_4, D_{24})$ &$13\cdot 7$ &3, 4, 12, 24 &1, 1, 3,
 2 \\ \hline
10 &$(PSL(2, 59), A_5, D_{60})$ &$59\cdot 29$ &6, 10, 15, 30, 60
 &1, 1, 1, 12, 22 \\ \hline
11 &$(PSL(2, 61), A_5, D_{60})$ &$61\cdot 31$ &6, 10, 15, 30, 60
 &1, 1, 1, 12, 25 \\ \hline
12 &$(PSL(2, 23), S_4, D_{24})$ &$23\cdot 11$ &3, 4, 6, 12, 24
 &1, 1, 1, 6, 7 \\ \hline
13 &$(M_{23}, PSL(3,4):Z_2, Z_2^4:A_7)$ &$23\cdot 11$ &21, 112, 120
 &1, 1, 1 \\ \hline
\end{tabular}}
\end{center}

\begin{center}
TABLE 2: Bi-coset graphs $\BB(G,L,R;$ $D)$

deriving non-biprimitive semisymmetric graphs
\vskip 5mm
{\footnotesize
 \begin{tabular}{|c|c|c|c|c|} \hline
Row &$(G,L,R)$ &$pq$ & Degree of $u$ & Number \\ \hline \hline
 1 &$(S_7,S_6,S_3\times S_4)$ &$7\cdot 5$ &15, 20 &1, 1\\ \hline
2 &$(S_p,S_{p-1},S_2\times S_{p-2})$ $(p\ge 5)$ &$p\cdot \frac {p-1}2$
 &$p-1, \frac {p^2-3p+2}2$ &1, 1 \\ \hline
3 &$(PSL(5, 2),Z_2^4:GL(4,2),Z_2^6:(S_3\times GL(3, 2))$ &$31\cdot 5$
&15, 140 &1, 1 \\ \hline
4 &$(PSL(2, 11),A_5,D_{12})$ &$11\cdot 5$
 &10, 15, 30 &1, 1, 1 \\ \hline
5 &$(PSL(2, 11),A_5,A_4)$ &$11\cdot 5$
 &20, 30 &1, 1 \\ \hline
6 &$(Z_p:Z_{qk},Z_{qk},Z_k) $ &$p\cdot q$ &$qk$ &1 \\
 &$qk\di (p-1), k\ge 2, q<\frac {p-1}2$ & & & \\ \hline
\end{tabular}}
\end{center}

\vskip 5mm
\noindent
{\it Note: in rows 4, 5 and 6 of TABLE 1, $\s $ is an involution
in $PGL(2,p)\setminus PSL(2,p)$, where $p=59$, 61 and 23 respectively,
and $\Aut \BB (G,L,R;D)\cong G = PSL(2,p)$.}

\section{Proof of Main  Results}
First we prove Theorem~\ref{a} and Corollary~\ref{b}.

\vskip 3mm

\demo
 Using {\sc Magma}\cite{BCP},   Conder and Poto\v{c}nik generate all the connected
semisymmetric cubic graphs of order up to 10000, see \cite{CP}.  With their  result,   we check that  every
semisymmetric cubic  graph of order up to 3000 has a  Hamilton cycle,
with the help of {\sc Magma}. So Theorem~\ref{a} is proved.

It was shown in \cite{II} that there are five cubic biprimitive semisymmetric graphs with respective order  110, 126, 182, 506 and 990, and automorphism group $PGL_2(11)$, $G_2(2)$, $PGL_2(13)$, $PSL_2(23)$ and $\mathrm{Aut}(M_{12})$. So by  Theorem~\ref{a}, every such graph has a  Hamilton cycle and thus
  Corollary~\ref{c} is proved.\qed

\vskip 3mm

 \vskip 3mm

To prove Theorem~\ref{c}, suppose that $Y$ is a semisymmetric graph of order $2pq$, where $q < p$, with bipartition $V(Y)=U(Y)\cup W(Y)$. Let $\A=\Aut(Y)$.
Then we shall show every graph in TABLEs 1 and 2 has a Hamilton cyclic.  We shall divide the proof into two cases: non-biprimitive cases and biprimitive cases, separately.

\subsection{Non-biprimitive cases}

Suppose $Y$ is not biprimitive. Then $Y\cong \ux$,  where $\ux$  is derived from the graphs  $X$   given in TABLE 2. First we find a   semiregular automorphism $\ba$ of order $p$ for $\ux$.

\begin{lem}\label{semiregular}
  Let $X=\BB(G,L,R;D)$  a bi-coset graph in  In TABLE 2.  Then   its  derived semisymmetric graph $\ux$ admits a semiregular automorphism $\ba$ of order $p$.
\end{lem}
\demo
Take an element $t$ of order $p$ in $G$. Since $p \nmid |L|$ and $p\nmid |R|$,  $\lg t\rg $ is transitive on $U(X)$ and acts semiregularly on $W(X)$. Then
$U(X)=\{Lt^j \di j\in \ZZ_p\}$ and $W(X)=[G:R]$ is a union of $q$  orbits of $\lg t\rg$.
Let $\Bc=\{B_i \di B_i $ , $i \in \ZZ_q\}$  be the set of  orbits of $\lg t\rg $  on $W(X)$, and pick up a vertex $Ra_i$ for each $B_i$ so that  $B_i=\{Ra_it^j\di j\in \ZZ_p\}$.

Now, the derived graph $\ux$ of $X$  is  described     as follows:
$$U(\ux)=\{(Lt^i,j)\di i \in \ZZ_q, j\in \ZZ_p \},~~~~ W(\ux)=W(X),$$
$$E(\ux)=\{\{Lt^i,j), Ra_lt^k\} \di \{Lt^j,   Ra_lt^k\} \in E(X),  i, k\in \ZZ_p, j, l\in \ZZ_q \}.$$
Define a bijection $\a$ by $t$ on $V(\ux)$ as follows:

$$(Lt^i,j)\to (Lt^{i+1},j),\quad  Ra_lt^k\to Ra_lt^{k+1},$$
 where $i, k\in \ZZ_p, j, l\in \ZZ_q$. Since $$\{(Lt^i,j), Ra_lt^k\}\in E(\ux)\Longleftrightarrow \{(Lt^{i+1},j), Ra_lt^{k+1}\}\in E(\ux),$$
 $\ba$ is an automorphism of graph $\ux$. Clearly, $\a $ acts semiregularly on $V(\ux)$.\qed

\begin{lem}\label{imp}
  Let $X$ and $\ux$ be the graphs as in last lemma.  Then $\ux$ has a Hamilton cycle.
\end{lem}
\demo   By  Lemma~\ref{semiregular}, $\lg \ba\rg$  has $q$ orbits on both $U(\ox)$ and $W(\ox)$, say ${\bf B}'=\{B'_i\di i\in \ZZ_p\}$ and   ${\bf B}=\{B_i\di i\in \ZZ_p\}$, respectively.
Let $\si$ be the block graph of $\ox$ induced by $\a$ with biparts ${\bf B}'$ and ${\bf B}$, where two orbits are adjacent if and only if there is at least  one edge in $\ux$
between two of them. Since  any vertex $Ra_lt^k$ is adjacent to a vertex in $U(X)$ in $X$,  $Ra_lt^k$ is adjacent to a vertex in each orbits in $U(X)$ in $X$. This implies $\si$ is
  a complete bipartite graph with $|V(\si)|=2q$. Obviously, $\si$ contains a Hamilton cycle. Pick up one such cycle: $B'_0B_0B'_1B_1\cdots B'_{q-1}B_{q-1}B'_0$.
   Then there exists a Hamilton path: $u_0w_0u_1w_1\cdots u_{p-1}w_{q-1}$ in $\ux$,  where $u_i\in B_i', w_i \in B_i$ for $i\in \ZZ_q$.
   Checking TABLE 2, one  may see that  $w_{q-1}$ is adjacent to at least two  vertices in $B_0'$. Since $\ba$ is transitive on $B_0$,  there exists a $d\in \lg \ba\rg $
   such that $\{d(u_0), w_{q-1})\} \in E(\ux)$ but $d$ moves $u_0$.
    For $k\in \ZZ_q$, define the following paths in $\ux$:
$$P_k=d^k(u_0), \d^k(w_0), d^k(u_1), d^k(w_1), \ldots, d^k(u_{q-1}), d^k(w_{q-1}).$$
Then $P_0P_1\cdots P_{q-1}u_0$ is a Hamilton cycle in $\ux$, as desired. \qed
\subsection{Biprimitive cases}
Suppose that $Y$ is biprimitive. Then $Y$ is one of bi-coset graphs given in TABLE 1. In what follows, $Y$ is said to be {\it bi-Cayley} if  $\Aut(Y)$ contains a subgroup acts regular on two biparts.  More precisely, we have the following two  cases essentially: $Y$ is a bi-Cayley graph or not.

\begin{lem} \label{bic} Suppose that $Y$ is a  bi-Cayley graph. Then it has a  Hamilton cyclic.
\end{lem}
\demo  Suppose that $Y$  a  bi-Cayley graph.  Checking TABLE 1, we know that $Y$ is one of graphs given in TABLE 3.
 Now $\Aut(Y)$ contains a subgroup $R$ acts regularly on two biparts. Since $|U(Y)|=|W(Y)|=pq$, where $p$ and $q$ are primes and $p>q$, we know that  $R$ is metacycic, that is
$$R=\lg a, b\di   a^p=b^q=1, a^b=a^r\rg ,$$
where $r^q\equiv 1(\mod p)$.
\begin{center}
TABLE 3: Bi-Cayley graphs

with vertex-transitive block graphs of order $2q$
\vskip 3mm
{\small
\begin{tabular}{|c|c|c|c|c|} \hline
Row &$(G, L, R)$ &$pq$ & Valency & Number \\ \hline \hline
1 &$(PSL(2, 59), A_5, A_5)$ &$59\cdot 29$ &60 &4 \\ \hline
3 &$(PSL(2, 23), S_4, S_4)$ &$23\cdot 11$ &12, 24 &1, 1 \\ \hline
4 &$(PSL(2, 59), A_5, A_5^{\s })$ &$59\cdot 29$ &60 &4 \\ \hline
6 &$(PSL(2, 23), S_4, S_4^{\s })$ &$23\cdot 11$ &12, 24 &1, 1 \\ \hline
7 &$(PGL(2, 11), S_4, D_{24})$ &$11\cdot 5$ &3, 4, 12, 24 &1, 1, 2, 1\\ \hline
10 &$(PSL(2, 59), A_5, D_{60})$ &$59\cdot 29$ &6, 10, 15, 30, 60
 &1, 1, 1, 12, 22 \\ \hline
12 &$(PSL(2, 23), S_4, D_{24})$ &$23\cdot 11$ &3, 4, 6, 12, 24
 &1, 1, 1, 6, 7 \\ \hline
13 &$(M_{23}, PSL(3,4):Z_2, Z_2^4:A_7)$ &$23\cdot 11$ &21, 112, 120
 &1, 1, 1 \\ \hline
\end{tabular}}
\end{center}

 Suppose $Y$ has valency $k$. Then $k\ge 3$. Take $u\in U(Y)$ and $w\in W(Y)$.   Suppose that $u$ is adjacent to
$\{w^{b^{i_l}a^{j_l}}\di l\in \ZZ_k\}$ and set $S=\{b^{i_l}a^{j_l}\di l\in \ZZ_k\}\subset R.$
Relabeling $w^{b^{i_0}a^{j_0}}$ by $w$, we may set $b^{i_0}a^{j_0}=1\in S$.  Noting that  $\lg a\rg $ is normal in $R$,  we consider   $R$-block graphs of $Y$ induced by $\lg a\rg $, which must be connected. Thus one of elements,
  say $b^{i_1}a^{j_1}\not\in \lg a\rg $  in $S$ is not contained in $\lg a\rg$. Then $|b^{i_1}a^{j_1}|=q$.
 Relabeling $b^{i_1}a^{j_1}$ by $b$, we set $b^{i_1}a^{j_1}=b$. In short, $1,b\in S$ and $u$ is adjacent to
$w$ and $w^b$. Take  into account, for any $g\in R$, the neighborhood of $u^g$ is $\{w^{sg}\di s\in S\}$ and the neighborhood of $w^g$ is $\{w^{s^{-1}g}\di s\in S\}$.
 We divide the proof into the following two cases, separately, depending on whether or not $u$ is adjacent to at least two vertices in a $R$-block, while the second case is call the {\it cover case}.

\vskip 3mm
{\it Case 1: Non-cover case}

(1) $S\cap \lg a\rg\ne \emptyset $:  Set $a^{k}\in S$, noting that $1, b\in S$.
Then
$$P_0=u,w^b,u^b,w^{b^2},\ldots, w^{b^{q-1}}, u^{b^{q-1}}, w$$ is a simple path.   For any $h\in \ZZ_p$,
$$P_h:=u^{a^{-kh}},w^{ba^{-kh}},u^{ba^{-kh}},w^{b^2a^{-kh}},\ldots, w^{b^{q-1}a^{-kh}}, u^{b^{q-1}a^{-kh}}, w^{a^{-kh}}$$ is a  simple cycle too. Clearly, for any $h_1\ne h_2$ in $\ZZ_p$, $P_{h_1}$ and $P_{h_2}$ have no common vertices.
 Then we get a Hamilton cyclic of $Y$ as follows:
$$P_0, P_1, P_2, \ldots , P_{p-1}, u.$$

\vskip 3mm
(2)  Suppose that $S\cap \lg a\rg=\emptyset $ but  $1, b^ja^i, b^ja^{i'}\in S$ where $i\ne j$: no loss, set $j=1$, $i'=0$ and $i\ne 0$. Set $x=ba^i$.
Observe a simple path:
$$P_0=w^x, u^x,w^{bx},u^{bx},w^{b^2x}, u^{b^2x},\ldots, w^{b^{q-1}x}, u^{b^{q-1}x},$$
noting that $b^{q-1}x=a^i$. For any $h\in \ZZ_p$,
$$P_h=w^{xa^{ih}}, u^{xa^{ih}},w^{bxa^{ih}},u^{bxa^{ih}},w^{b^2xa^{ih}}, u^{b^2xa^{ih}},\ldots, w^{b^{q-1}xa^{ih}}, u^{b^{q-1}xa^{ih}}$$
is a  simple cycle too. Then we get a Hamilton cyclic of $Y$ as follows:
$$P_0, P_1, P_2, \ldots , P_{p-1}, w^x,$$
noting that the ended vertex of $P_{p-1}$ is exactly  $u^{b^{q-1}xa^{i(p-1)}}=u$.

\vskip 3mm
{\it Case 2: Cover case}

Suppose that  $u$ is adjacent to at most  one  vertices in every  $R$-block.
Then  the valency of $Y$  is smaller   than $q+1$. Checking TABLE 3, we only need to consider the following cases (some of them might be noncover case):

\begin{center}
{\small
\begin{tabular}{|c|c|c|c|c|} \hline
Row &$(G, L, R)$ &$pq$ & Valency & Number \\ \hline \hline
7 &$(PGL(2, 11), S_4, D_{24})$ &$11\cdot 5$ &3, 4 &1, 1\\ \hline
10 &$(PSL(2, 59), A_5, D_{60})$ &$59\cdot 29$ &6, 10 , 15 &1, 1, 1 \\ \hline
12 &$(PSL(2, 23), S_4, D_{24})$ &$23\cdot 11$ &3, 4, 6
 &1, 1, 1\\ \hline
\end{tabular}}
\end{center}

For these cases, we find a Hamilton  cycle by Magma.
\qed

\begin{lem} \label{bic2} Suppose that $Y$ is a biprimitive bi-coset graph which is not bi-Cayley. Then it has a  Hamilton cyclic.
\end{lem}
\demo Checking TABLE 1, we know that $Y$ is given in following TABLE 4. For these cases,  we find a Hamilton  cycle by Magma. \qed

\begin{center}
TABLE 4: Biprimitive bi-coset graphs which are not bi-Cayley
\vskip 5mm
{\small
\begin{tabular}{|c|c|c|c|c|} \hline
Row &$(G, L, R)$ &$pq$ & Valency & Number \\ \hline \hline
2 &$(PSL(2, 61), A_5, A_5)$ &$61\cdot 31$ &60 &4 \\ \hline
5 &$(PSL(2, 61), A_5, A_5^{\s })$ &$61\cdot 31$ &60 &5 \\ \hline
8 &$(PSL(2, 13), A_4, D_{12})$ &$13\cdot 7$ &6, 12 &1, 2 \\ \hline
9 &$(PGL(2, 13), S_4, D_{24})$ &$13\cdot 7$ &24 & 2 \\ \hline
11 &$(PSL(2, 61), A_5, D_{60})$ &$61\cdot 31$ &6, 10, 15, 30, 60
 &1, 1, 1, 12, 25 \\ \hline
\end{tabular}}
\end{center}

From Lemmas \ref{bic} and \ref{bic2}, and with the help of {\sc Magma}\cite{BCP}, the following Lemma holds.
\begin{lem} \label{bip}
The semisymmetric graphs in TABLE 1, are hamiltonian.
\end{lem}

Finally, by Lemma \ref{imp} and \ref{bip}, one can get Theorem~\ref{c}.

\end{document}